\documentclass[11pt]{article}
\usepackage{}
\usepackage[total={6in, 9in}]{geometry}
\usepackage{amsfonts}
\usepackage{amsthm}
\usepackage{amssymb}
\usepackage{amsmath,amsfonts}
\usepackage{mathrsfs}
\usepackage{cases}
\usepackage{latexsym,bm}
\usepackage{indentfirst}
\usepackage{color}
\usepackage{ifpdf}
\usepackage{graphicx}
\usepackage{psfrag}
\usepackage{dsfont}
\usepackage[
pdfauthor={CESS},
pdftitle={Spanning trails},
pdfstartview=XYZ,
bookmarks=true,
colorlinks=true,
linkcolor=blue,
urlcolor=blue,
citecolor=blue,
bookmarks=true,
linktocpage=true,
hyperindex=true
]{hyperref}

\newtheorem{THM}{\textbf{Theorem}}

\newtheorem{LEM}{\textbf{Lemma}}[section]
\newtheorem{CLA}{\textbf{Claim}}[section]
\newtheorem{CON}{\textbf{Conjecture}}

\newtheorem{ASP}[CLA]{\textbf{Assumption}}

\newcommand{\pf}{\noindent\textbf{Proof}.\quad}

\newcommand{\qqed}{\hfill $\blacksquare$\vspace{1mm}}

\newcommand{\iC}{\overset{\leftharpoonup }{C}}
\newcommand{\iP}{\overset{\leftharpoonup }{P}}
\newcommand{\iD}{\overset{\leftharpoonup }{D}}
\newcommand{\iQ}{\overset{\leftharpoonup }{Q}}
\newcommand{\iR}{\overset{\leftharpoonup}{R}}
\newcommand{\oC}{\overset{\rightharpoonup }{C}}
\newcommand{\oP}{\overset{\rightharpoonup }{P}}
\newcommand{\oD}{\overset{\rightharpoonup }{D}}
\newcommand{\oQ}{\overset{\rightharpoonup }{Q}}
\newcommand{\oR}{\overset{\rightharpoonup }{R}}

\newcommand{\sF}{\mathcal{F}}
\newcommand{\Nb}{\overline{V}}

\begin{document}
\title{Hamiltonian cycles  in 3-tough $2K_2$-free graphs}
\author{ Songling Shan\thanks {Department of Mathematics, Vanderbilt University, Nashville, TN 37240,
U.S.A. ({\tt songling.shan@vanderbilt.edu}).}
}

\date{\today}
\maketitle

\emph{\textbf{Abstract}.}
A graph is called $2K_2$-free if it does not contain two independent edges as an induced subgraph.
Broersma, Patel, and Pyatkin showed that every 25-tough $2K_2$-free
graph with at least three vertices is hamiltonian.  In this paper,
we improve the required toughness in this result from 25 to 3.
%
%

\emph{\textbf{Keywords}.} Toughness; Hamiltonian cycle; $2K_2$-free graph

\vspace{2mm}

\section{Introduction}

Graphs considered in this paper are simple, undirected,  and finite.
Let $G$ be a graph.
Denote by $V(G)$ and  $E(G)$ the vertex set and edge set of $G$,
respectively. For $v\in V(G)$,  $N_G(v)$ denotes the set of neighbors
of $v$ in $G$. For $S\subseteq V(G)$,
$N_G(S)=\bigcup_{x\in S}N_G(x)-S$.
For $H\subseteq G$  and $x\in V(G)$,
define $V_H(x)=N_G(x)\cap V(H)$ and $V_H(S)=N_G(S)\cap V(H)$.
Let $S\subseteq V(G)$.  Then the subgraph induced by $V(G)-S$ is denoted by
$G-S$. For notational simplicity we write $G-x$ for $G-\{x\}$.
If $uv\in E(G)$ is an edge, we write $u\sim v$.
Let $V_1,
V_2\subseteq V(G)$ be two disjoint vertex sets. Then $E_G(V_1,V_2)$ is the set
of edges of $G$ with one end in $V_1$ and the other end in $V_2$.

The number of components of $G$ is denoted by $c(G)$. Let $t\ge 0$ be a
real number. The graph is said to be {\it $t$-tough\/} if $|S|\ge t\cdot
c(G-S)$ for each $S\subseteq V(G)$ with $c(G-S)\ge 2$. The {\it
toughness $\tau(G)$\/} is the largest real number $t$ for which $G$ is
$t$-tough, or is  $\infty$ if $G$ is complete. This concept, a
measure of graph connectivity and ``resilience'' under removal of
vertices, was introduced by Chv\'atal~\cite{chvatal-tough-c} in 1973.
 It is  easy to see that  if $G$ has a hamiltonian cycle
then $G$ is 1-tough. Conversely,
 Chv\'atal~\cite{chvatal-tough-c}
 conjectured that
there exists a constant $t_0$ such that every
$t_0$-tough graph is hamiltonian.
 Bauer, Broersma and Veldman~\cite{Tough-counterE} have constructed
$t$-tough graphs that are not hamiltonian for all $t < \frac{9}{4}$, so
$t_0$ must be at least $\frac{9}{4}$.

There are a number of papers on
 Chv\'atal's toughness conjecture,
and it  has
been verified when restricted to a number of graph
classes~\cite{Bauer2006},
including planar graphs, claw-free graphs, co-comparability graphs, and
chordal graphs. A
graph $G$ is called {\it $2K_2$-free\/} if it does not contain two independent
edges as an induced subgraph.
Recently, Broersma, Patel and Pyatkin~\cite{2k2-tough} proved that
every 25-tough $2K_2$-free graph on at least three vertices is
hamiltonian.


%


The class of $2K_2$-free graphs is well studied, for instance, see~\cite{
2k2-tough, CHUNG1990129, MR845138, 2k21, 1412.0514, MR2279069,  MR1172684}.
It is a superclass of {\it split\/} graphs,
where the vertices can be partitioned into a clique and an independent set.
One can also easily check that every {\it cochordal\/} graph (i.e., a graph that is the complement of a
chordal graph) is $2K_2$-free and so the class of $2K_2$-free graphs is
at least as rich as the class of chordal graphs.
In~\cite{2k21}, Gao and Pasechnik proposed the following  conjecture.


\begin{CON}\label{2k2h}
Every $2$-tough $2K_2$-free graph with at least three vertices is hamiltonian.
\end{CON}

In this paper, we support  Conjecture~\ref{2k2h} as well as improve the main result in~\cite{2k2-tough} by showing the following result.

\begin{THM}\label{main}
Let $G$ be a $3$-tough $2K_2$-free graph with at least three vertices. Then $G$
is hamiltonian.
\end{THM}

In~\cite{MR1392734} it was shown that every 3/2-tough split graph on at least three vertices
is hamiltonian.  And the authors constructed
a sequence $\{G_n\}_{n=1}^{\infty}$ of split graphs with no 2-factor and $\tau(G_n)\rightarrow 3/2$.
So $3/2$ is the best possible toughness for split graphs to be hamiltonian.
Since split graphs are $2K_2$-free,
we cannot decrease the bound in Theorem~\ref{main} below 3/2.
Although we are not sure about the best possible toughness for guaranteeing  $2K_2$-free
graphs to be hamiltonian, we believe that Conjecture~\ref{2k2h}
might be true. In fact, in the proof of Theorem~\ref{main},
except for one case where 3-tough is needed, all other cases only need
a toughness of 2.



\section{Proof of Theorem~\ref{main}}


We need the following lemma for the existence of a 2-factor in a graph.

\begin{LEM}[Enomoto et al.~\cite{MR785651}]\label{2-factor}
Every $k$-tough graph has a $k$-factor
if $k|V(G)|$ is even and $|V(G)|\ge k+1$.
\end{LEM}

We will also need some notation. Let $C$ be an oriented cycle. For $x\in V(C)$,
denote the successor of $x$ by $x^+$ and the predecessor of $x$ by $x^-$.
Let $S\subseteq V(C)$ be an independent set in $C$. Then $S^+=\{x^+\,|\, x\in S\}$,
and $S^-$ is defined similarly. Let $D$ be another oriented cycle disjoint with $C$
and $T\subseteq V(D)$ be an independent set in $D$. Then $(S\cup T)^+=S^+\cup T^+$
and $(S\cup T)^-=S^-\cup T^-$.
For $u,v\in V(C)$, $u\oC v$ denotes the portion of $C$
starting at $u$, following $C$ in the orientation,  and ending at $v$.
Likewise, $u\iC v$ is the opposite portion of $C$ with endpoints as $u$
and $v$. Given two vertex-disjoint cycles $C$ and $D$. Suppose  $P_c$ is  a portion of
$C$ with endpoints $u,v$ and $P_d$ is a portion of $D$ with endpoints $x,y$.
If $v$ and $x$ are adjacent, we write $uP_cvxP_dy$ as
the concatenation of $P_c$ and $P_d$ through the edge $vx$.  We will assume all cycles
in consideration are oriented.


\proof[Proof of Theorem~\ref{main}]
The graph $G$ is 3-tough, so it has a 2-factor by
Lemma~\ref{2-factor}. We take a 2-factor of $G$ such that
it contains as few cycles as possible.
Let $\sF$ be the set of cycles in such a 2-factor. We may assume that
$\sF$ contains at least two cycles. For otherwise, the only cycle in
$\sF$ is a hamiltonian cycle of $G$.

Let $x\in V(G)$ be a vertex. As cycles in $\sF$
form a 2-factor of $G$, there exists a unique cycle, say
$C\in \sF$ such that
 $x\in V(C)$.
If there exists a   cycle
$D\in \sF-\{C\}$ such that $x$ is adjacent to
two consecutive vertices on $D$ in $G$,
 we say  $x$ is of  {\it A-type}\,(w.r.t. $D$). If $x$
 is not of A-type w.r.t. any cycles in $\sF-\{C\}$, we say
 $x$ is of {\it B-type}. Denote
$$
   A=\{x\in V(G)\,|\, \mbox{$x$ is of A-type}\}\quad\quad \mbox{and}\quad\quad  B=V(G)-A.
 $$
Let  $xy\in E(C)$ be an edge.  We say $xy$ is of {\it A-type} if $x,y\in A\cap V(C)$;
we say $xy$ is of {\it B-type} if $x,y\in B\cap V(C)$; otherwise,
$xy$ is of {\it AB-type}. We say $C$ is {\it AB-alternating} if all edges of $C$
are of AB-type. It is clear that if $C$ is AB-alternating, then $C$ is an even cycle.
For a cycle $D\in \sF-\{C\}$ and the edge $xy\in E(C)$, we denote
$$
V_D(xy)=V_D(x)\cup V_D(y)\quad\quad \Nb_D(xy)=V(D)-V_D(xy),
$$
where recall that $V_D(x)=N_G(x)\cap V(D)$. 

\begin{CLA}\label{nonadjacency}
Let $C,D\in \sF$ be two distinct cycles. If $x\in V(C)$
has a neighbor $u\in V(D)$,  then $x^+\not\sim u^+, u^-$ and $x^-\not\sim u^+, u^-$.
\end{CLA}

\pf Assume on the contrary that $x^+\sim u^+$. Then $xu\iD u^+x^+\oC x$
combines $C$ and $D$ into a single cycle. This gives a contradiction to
the minimality of $|\sF|$. Similar construction shows that $x^+\not\sim u^-$ and
$x^-\not\sim u^+, u^-$.
\qed

\begin{CLA}\label{no-a-type-edge}
 No cycle in $\sF$ containing an A-type edge.
\end{CLA}
\pf Assume on the contrary that there exists $C\in\sF$
and $xy\in E(C)$ such that $xy$ is of A-type.  Let $D,Q\in \sF-\{C\}$
such that $x\sim u, u^+$ with $uu^+\in E(D)$ and $y\sim v, v^+$
with $vv^+\in E(Q)$.
As $x\sim u, u^+$, $y\not\sim u, u^+$ by Claim~\ref{nonadjacency}.
 Let $z$ be the other neighbor of $y$ on $C$.
 Then $z\sim u$ or $z\sim u^+$ by considering
 the two independent edges $yz$ and $uu^+$. By reversing the
 orientation of $C, D$ if necessary, we assume that $y=x^+$ and
 $z\sim u$.
  Then
 $$
 \left\{
   \begin{array}{ll}
     xu^+\oD v yv^+\oD uz \oC x \quad \mbox{combines $C$ and $D$ into one cycle}, & \hbox{if $D=Q$;} \\
     xu^+D uz \oC x \quad \mbox{and} \quad vyv^+\oQ v  \quad \mbox{integrate  $C, D, Q$ into two cycles}, & \hbox{if $D\ne Q$.}
   \end{array}
 \right.
$$
\qed

 \begin{CLA}\label{I-xy-alter-size-inde}
Let $C\in \sF$,  $xy\in E(C)$.
Denote
$$
I_{xy}=\bigcup_{D\in \sF-\{C\}} \Nb_D(xy).
$$
Then each of the following holds.
\begin{enumerate}
\item [$(1)$] $I_{xy}$ is an independent set in $G$.
  \item [$(2)$] If $xy$ is of B-type, then for any $D\in \sF-\{C\}$, vertices on $D$ are alternating between
  $I_{xy}$ and $V(G)-V(C)-I_{xy}$.
  \item [$(3)$] If $xy$ is of B-type, then $|I_{xy}|=\frac{1}{2}|V(G)-V(C)|$.

\end{enumerate}
\end{CLA}
\pf  To show $(1)$, assume on the contrary that
there exist $u,v\in I_{xy}$ such that
$u\sim v$. Then $E_G(\{x,y\},\{u,v\}) \ne \emptyset$
by the $2K_2$-freeness assumption of $G$.
Consequently, at least one of $u$ and $v$ is not an element in $I_{xy}$.
This gives a contradiction.

Assume now that $xy$ is a B-type edge.
Let $D\in \sF-\{C\}$. We show that for any edge
$uv\in E(D)$, there is one and exactly one vertex in $\{u,v\}$
is in $V_D(xy)$. One of $u,v$ must be in $V_D(xy)$ is
guaranteed by the $2K_2$-freeness of $G$.
Suppose, w.l.o.g., that $u\in V_D(xy)$ with $u\sim x$.
Then by Claim~\ref{nonadjacency}, $v\not\sim y$. As $x$ is of B-type,
we further know that $v\not\sim x$. Thus, $v\in \Nb_D(xy)$.  This gives $(2)$.
The statement $(3)$ is an immediate consequence of $(2)$.
\qed

%

\begin{CLA}\label{a-plus-ind}
Let $A^+$ be the set of successors of vertices in $A$.
Then $A^+$ is an independent set in $G$.
\end{CLA}

\pf Suppose on the contrary that there exist $x^+, y^+\in A^+$
with $x^+y^+\in E(G)$.

Assume $x^+\in V(C)$ with predecessor as $x$, $y^+\in V(D)$
with predecessor as $y$ for cycles $C, D\in \sF$.  Then both $x$ and $y$ are A-type vertices.
Let $Q, R\in \sF$
with $uu^+\in E(Q)$ and $vv^+\in E(R)$
such that $x\sim u, u^+$ and $y\sim v, v^+$.
As $x\sim u, u^+$ and $y\sim v, v^+$, we know that $x^+\not\sim u^-, u, u^+, u^{++}$ and $y^+\not\sim v^-, v, v^+, v^{++}$
by Claim~\ref{nonadjacency}. Since $x^+y^+\in E(G)$, by the $2K_2$-freeness of $G$,
$y^+$ is adjacent to one of $u,u^+$ and $x^+$ is adjacent to one of $v, v^+$.
Thus,  $\{u,u^+\}\cap \{v, v^+\}=\emptyset$.
We consider two cases
for completing the proof.

 \smallskip\noindent
{\bf Case \ref{a-plus-ind}.1:} $C=D$.

\smallskip\noindent
{\bf Case \ref{a-plus-ind}.1.1:} $C=D$ and $Q=R$.

We combine $C$ and $Q$ into a single cycle as follows.
$$
\left\{
  \begin{array}{ll}
    xu\iQ v^+y \iC x^+ v \iQ u^+y^+\oC x, & \hbox{if $x^+\sim v, y^+\sim u^+$;} \\
   xu^+\oQ vx^+ \oC y v^+ \oQ uy^+\oC x, & \hbox{if $x^+\sim v, y^+\sim u$;} \\
    xu\iQ v^+x^+ \oC y v \iQ u^+y^+\oC x, & \hbox{if $x^+\sim v^+, y^+\sim u^+$;} \\
   xu^+\oQ vy \iC x^+ v^+ \oQ uy^+\oC x, & \hbox{if $x^+\sim v^+, y^+\sim u$.}
  \end{array}
\right.
$$

\smallskip\noindent
{\bf Case \ref{a-plus-ind}.1.2:} $C=D$ and $Q\ne R$.

Recall that $\{u, u^+\}\cap \{v, v^+\}=\emptyset$. Thus, $E_G(\{u, u^+\}, \{v, v^+\})\ne \emptyset$.
By reversing the orientations of $Q$ and $R$
if necessary, we assume that $u\sim v$.
Then
$$
xu^+\oQ uv \iR v^+y \iC x^+y^+\oC x
$$
combines $C$,  $Q$ and $R$ into a single cycle.
%

\smallskip\noindent
{\bf Case \ref{a-plus-ind}.2:} $C\ne D$.

\smallskip\noindent
{\bf Case \ref{a-plus-ind}.2.1:} $C\ne D$ and $Q=R$.

As $Q\ne C$ and $R\ne D$ by the definition of A-type vertices, we have $Q\not\in \{C,D\}$.
 Recall that $\{u,u^+\}\cap \{v, v^+\}=\emptyset$.
 Thus, $E_G(\{u, u^+\},\{v,v^+\})\ne \emptyset$,
by the $2K_2$-freeness of $G$.
By reversing the orientation of $Q$ if necessary, we assume $y^+\sim u$.
Then $uv^+\not\in E(Q)$. As otherwise, $u=v^{++}$ and so $y^+\sim v^{++}$.
However, $y^+\not\sim v^{++}$ by the argument prior to Case \ref{a-plus-ind}.1.

We cover vertices in $V(C)\cup V(D)\cup V(Q)$ by one or two
cycles as below.
$$
\left\{
  \begin{array}{ll}
    xu^+\iQ v y\iD y^+ x^+  \oC x, & \hbox{if $u^+v\in E(Q)$;} \\
    xu\iQ v^+ y\iD y^+ x^+  \oC x, u^+\oQ vu^+, & \hbox{if $u^+\sim v$ but $u^+v\not\in E(Q)$;} \\
    xu \iQ v^+ u^+ \oQ  vy \iD y^+ x^+ \oC x, & \hbox{if $u^+\sim v^+$;} \\
    xu^+ \oQ v u \iQ  v^+y \iD y^+ x^+ \oC x, & \hbox{if $u\sim v$;} \\
     xu^+ \oQ v y \iD y^+ x^+ \oC x,  u\iQ v^+ u,  & \hbox{if $u\sim v^+$.}
  \end{array}
\right.
$$

\smallskip\noindent
{\bf Case \ref{a-plus-ind}.2.2:} $C\ne D$ and $Q\ne R$.

As $x^+\not\sim u, u^+$, and $x^+\sim y^+$,
we get $y^+\not\in \{u, u^+\}$.  Consequently, $y\not\in\{u^-, u\}$.
Similarly,
$x^+\not\in \{v, v^+\}$ and  $x\not\in\{v^-, v\}$.

\smallskip\noindent
{\bf Case \ref{a-plus-ind}.2.2.1:} $C\ne D$,  and $Q=D, R=C$.

Again $E_G(\{u, u^+\},\{v,v^+\})\ne \emptyset$ and $E_G(\{u^+, u^{++}\},\{v,v^+\})\ne \emptyset$
by the $2K_2$-freeness of $G$.
We combine $C$ and $D$ into a single cycle as follows.
$$
\left\{
  \begin{array}{ll}
       xu\iD y^+ x^+  \oC vu^+\oD y v^+\oC x, & \hbox{if $u^+\sim v$;} \\
    xu\iD y^+ x^+  \oC vy\iD u^+ v^+\oC x, & \hbox{if $u^+\sim v^+$;} \\
    xu^+\oD yv  \iC x^+y^+\iD u v^+\oC x, & \hbox{if $u\sim v^+$;} \\
     xu^+ \iD y^+ x^+ \oC v u^{++}\iD y v^+ \oC x,  & \hbox{if $u\sim v,\, u^+\not\sim v, v^+$, and $ u^{++}\sim v$;}\\
      xu^+ \iD y^+ x^+ \oC v y\iD u^{++} v^+ \oC x,  & \hbox{if $u\sim v,\, u^+\not\sim v, v^+$, and $ u^{++}\sim v^+$.}
  \end{array}
\right.
$$

\smallskip\noindent
{\bf Case \ref{a-plus-ind}.2.2.2:} $C\ne D, Q\ne R$,  and $Q=D, R\not\in \{C,D\}$.

We cover vertices in $V(C)\cup V(D)\cup V(R)$ by one or two
cycles as below.
$$
\left\{
  \begin{array}{ll}
    xu^+\oD yv^+  \oR vu\iD y^+x^+\oC x, & \hbox{if $u\sim v$;} \\
    xu^+\oD yv  \iR v^+u\iD y^+x^+\oC x, & \hbox{if $u\sim v^+$;} \\
    xu\iD  y^+ x^+ \oC x, u^+ \oD y v^+ \oR vu^+ & \hbox{if $u^+\sim v$;} \\
      xu\iD  y^+ x^+ \oC x, u^+ \oD y v \iR v^+u^+ ,  & \hbox{if $u^+\sim v^+$.}
  \end{array}
\right.
$$
\smallskip\noindent
{\bf Case \ref{a-plus-ind}.2.2.3:} $C\ne D, Q\ne R$,  and $R=C, Q\not\in \{C,D\}$.

This case is symmetric to Case~\ref{a-plus-ind}.2.2.2, so we skip its proof.

\smallskip\noindent
{\bf Case \ref{a-plus-ind}.2.2.4:} $C\ne D, Q\ne R$,  and $Q\ne D, R\ne C$.

By reversing the orientations of $Q$ and $R$ if necessary, we assume that $u\sim v$.
Then
$$
xu^+\oQ uv  \iR v^+y\iD y^+x^+\oC x
$$
combines $C, D, Q, R$ into a single cycle.
\qed

\begin{CLA}\label{one-cycle-have-B-edge}
 We may assume that $\sF$ contains exactly one cycle $C$
such that $C$ has a B-type edge, and all other cycles in $\sF$ are AB-alternating.
\end{CLA}

\pf By Claim~\ref{a-plus-ind} that $A^+$ is an independent set in $G$, we know
that not all cycles in $\sF$ are AB-alternating. As otherwise,
let $S=A$. Then $c(G-S)=|A^+|=|A|=|S|$. We get that
$\tau(G)\le \frac{|S|}{|c(G-S)|}=1<3$. This gives a contradiction.

We then claim that  $\sF$ contains no
two cycles, say $C$ and $D$ both containing  a B-type edge. Assume on the contrary that
both $C$ and
$D$ contain a B-type edge. Suppose, w.l.o.g., that $|V(C)|\le |V(D)|$.
Let $xy\in E(C)$ be of B-type. By Claim~\ref{I-xy-alter-size-inde}, $I_{xy}$, the set of
non-neighbors of $x$ and $y$ in $V(G)-V(C)$, is an independent
set in $G$, Thus, $I_{xy}\cup \{x\}$ is also an independent set in $G$.
Let $S=V(G)-(I_{xy}\cup \{x\})$. Then $G-S$
has $|I_{xy}\cup \{x\}|$ components, each being an isolated vertex.
Note that  $|I_{xy}|=\frac{|V(G)-V(C)|}{2}=\frac{|V(G)-V(C)-V(D)|}{2}+\frac{|V(D)|}{2}$ by Claim~\ref{I-xy-alter-size-inde},
and $|V(C)|\le |V(D)|$.
Thus,
\begin{eqnarray*}
  \tau(G) &\le &  \frac{|S|}{c(G-S)} \\
   &=& \frac{|V(C)|-1+\frac{|V(G)-V(C)|}{2}}{|I_{xy}|+1} \\
   &=&  \frac{|V(C)|-1+\frac{|V(G)-V(C)-V(D)|}{2}+\frac{|V(D)|}{2}}{\frac{|V(G)-V(C)-V(D)|}{2}+\frac{|V(D)|}{2}+1}\\
   &\le & \frac{|V(D)|-1+\frac{|V(G)-V(C)-V(D)|}{2}+\frac{|V(D)|}{2}}{\frac{|V(G)-V(C)-V(D)|}{2}+\frac{|V(D)|}{2}+1}\\
   &= & \frac{\frac{|V(G)-V(C)-V(D)|}{2}+\frac{3|V(D)|}{2}-1}{\frac{|V(G)-V(C)-V(D)|}{2}+\frac{|V(D)|}{2}+1}\\
   & <& 3,
\end{eqnarray*}
showing a contradiction to the  assumption that $\tau(G)\ge 3$. (In fact, this is the only case where  3-tough is used.)
\qed

\begin{ASP}\label{C-only-B-edge}
We now fix $C\in \sF$ to denote the cycle which contains
a B-type edge, and assume that all other cycles in $\sF-\{C\}$
are AB-alternating.
\end{ASP}

\begin{CLA}\label{full-b-neighbors}
Let $D\in \sF-\{C\}$ and $xy\in E(C)$ be of B-type.
Assume that $V_D(xy)\cap B\ne\emptyset$.
Then $V_D(xy)=B\cap V(D)$, and either $V_D(x)=\emptyset$ and $V_D(y)=B\cap V(D)$ or
$V_D(x)=B\cap V(D)$ and $V_D(y)=\emptyset$.
\end{CLA}
\pf
Recall that $\Nb_D(xy)$ is an independent set in $G$, vertices on $D$
are alternating between $\Nb_D(xy)$ and $V_D(xy)$ by $(2)$
of Claim~\ref{I-xy-alter-size-inde}. Because $D$ is AB-alternating,
we then get $V_D(xy)=B\cap V(D)$ if $V_D(xy)\cap B\ne\emptyset$.
And so if $V_D(x)=\emptyset$, then $V_D(y)=B\cap V(D)$; and
if $V_D(y)=\emptyset$, then $V_D(x)=B\cap V(D)$. Thus, we
only show that either $V_D(x)$ or $V_D(y)$ has to be empty.

Assume to the contrary that $V_D(x)\ne \emptyset$ and $V_D(y)\ne\emptyset$.
As vertices on $D$
are alternating between $\Nb_D(xy)$ and $V_D(xy)$,
we can choose $u\in V_D(x)$ so that $u^{++}\in V_D(y)$.
Then $u^+\in A\cap V(D)$. Assume that
$u^+$ is of A-type w.r.t. $Q\in \sF-\{D\}$, i.e.,
$u^+\sim v, v^+$ with $vv^+\in E(Q)$.  Assume, w.l.o.g., that $y=x^+$.
As $x\sim u$ and $y\sim u^{++}$, we have that $u^+\not\sim x,y$ by Claim~\ref{nonadjacency}.
Thus, $\{v, v^+\}\cap \{x,y\}=\emptyset$.

\smallskip
 \noindent
{\bf  Case~\ref{full-b-neighbors}.1:} $Q=C$.

We combine $C, D$ into a single cycle as $xu \iD u^{++}y \oC v u^+ v^+ \oC x$.

\smallskip
 \noindent
{\bf  Case~\ref{full-b-neighbors}.2:}  $Q\ne C$.

We cover $V(C)\cup V(D)\cup V(Q)$ by two cycles
 as $xu \iD u^{++}y \oC x$ and $v \oQ v^+ u^+ v$.
\qed

\begin{CLA}\label{full-b-neighbors2}
Let $ D\in \sF-\{C\}$ and $x\in V(C)$.
Assume that $V_D(x)=B\cap V(D)$, then $V_D(x^+)=V_D(x^-)=\emptyset$.
\end{CLA}

\pf Note first that $N_G(x^+)\cap (B\cap V(D))=\emptyset$
and $N_G(x^-)\cap (B\cap V(D))=\emptyset$. As otherwise,
some vertex in $B\cap V(D)$ is adjacent to both vertices
in $\{x, x^+\}$ or $\{x,x^-\}$. This implies that
the vertex is of A-type w.r.t. $C$. Then we observe
that neither $x^+$ nor $x^-$
is adjacent to any vertex in $A\cap V(D)$ by Claim~\ref{nonadjacency}.
\qed

\begin{CLA}\label{v-0-ind-b}
Let  $x\in V(C)$.
Assume there exists $D\in \sF$
so that  $V_D(x)=B\cap V(D)$.  Then
$\{x^+\}\cup A^+$  is an independent set in $G$.
\end{CLA}

\pf  
As $A^+$ is already an independent set in $G$ by Claim~\ref{a-plus-ind},
we assume on the contrary that there exists $w\in A^+$
so that $x^+\sim w$.  Note that $x\ne w$, since $V_D(x)=B\cap V(D)$
and $V_D(w)\cap B=\emptyset$.
Assume $w\in V(Q)$ for some cycle $Q\in \sF$.
 Then the predecessor  $w^-$ of $w$
on $Q$ is of A-type.  Note that $V_D(x^+)=\emptyset$ by Claim~\ref{full-b-neighbors2} and
$x^+\sim w$ implies that $Q\ne D$.
Let $R\in \sF-\{Q\}$ with $vv^+\in E(R)$
so that $w^-\sim v, v^+$.
Let $z\in A\cap V(D)$.
As $V_D(x^+)=\emptyset$ and $x^+\sim w$, $w$
is adjacent to one of $z$ and $z^+$ by the $2K_2$-freeness of $G$.
Since $D$ is AB-alternating by Assumption~\ref{C-only-B-edge}
and $z\in A\cap V(D)$, $z^+\in B\cap V(D)$.
We see that $w\sim z$, because both $w, z^+\in A^+$
and $A^+$ is an independent set in $G$ by Claim~\ref{a-plus-ind}.
As $w^-\sim v, v^+$, $w\not\sim v, v^+$ by Claim~\ref{nonadjacency}.
Thus, $E_G(\{w^+\}, \{v, v^+\})\ne \emptyset$.
We consider two cases for finishing the proof.

\smallskip
 \noindent
{\bf  Case~\ref{v-0-ind-b}.1:} $Q\ne C$.

As $x^+\sim w$, we have $w^-\not\sim x^{++}, x$ by Claim~\ref{nonadjacency}.
Since $w^-\sim v, v^+$,  we then have that $v, v^+\not\in \{x, x^+, x^{++}\}$.

\smallskip
 \noindent
{\bf  Case~\ref{v-0-ind-b}.1.1:} $Q\ne C$ and $R=C$.

We combine $C,D, Q$ into one single cycle as below.
$$
\left\{
  \begin{array}{ll}
    x^+wz\oD z^-x \iC v^+ w^- \iQ w^+ v \iC x^+ , & \hbox{if $w^+\sim v$;} \\
     x^+wz\oD z^-x \iC v^+ w^+\oQ w^- v \iC x^+, & \hbox{if $w^+\sim v^+$.}
  \end{array}
\right.
$$
\smallskip
 \noindent
{\bf  Case~\ref{v-0-ind-b}.1.2:} $Q\ne C$ and $R= D$.

By the assumption, $V_D(x^+)=\emptyset$; particulary, $x^+\not\sim v, v^+$.
Since $w^-\sim v, v^+$, $w\not\sim v, v^+$ by Claim~\ref{nonadjacency}.
But $x^+w$ and $vv^+$ are two induced
disjoint edges.  This gives a contradiction to
the $2K_2$-freeness.

%
\smallskip
 \noindent
{\bf  Case~\ref{v-0-ind-b}.1.2:} $Q\ne C$ and $R\not\in \{C, D\}$.

Since $w^-\sim v, v^+$, $w\not\sim v, v^+$ by Claim~\ref{nonadjacency}.
Thus, $x^+\sim v$ or $x^+\sim v^+$.
By reversing the orientation of $R$ if necessary, we assume $x^+\sim v$.
Since $V_D(x)=B\cap V(D)$ and $z^+\in B\cap V(D)$, $x\sim z^+$.
Then
$$
xz^+\oD z w\oQ w^-v^+\oR vx^+ \oC x
$$
is a cycle which contains all the vertices  in $V(C)\cup V(D)\cup V(Q)\cup V(R)$.

\smallskip
 \noindent
{\bf  Case~\ref{v-0-ind-b}.2:} $Q=C$.

As $z\sim w$, $w^- \not\sim z^+, z^-$ by Claim~\ref{nonadjacency}.
Since $w^-\sim v, v^+$, we then get that
$v, v^+\not\in \{z^-, z, z^+\}$.

\smallskip
 \noindent
{\bf  Case~\ref{v-0-ind-b}.2.1:} $Q=C$ and $R=D$.

By the assumption, $V_D(x^+)=\emptyset$; particulary, $x^+\not\sim v, v^+$.
Since $w^-\sim v, v^+$, $w\not\sim v, v^+$ by Claim~\ref{nonadjacency}.
But $x^+w$ and $vv^+$ are two induced
disjoint edges.
This gives a contradiction to
the $2K_2$-freeness.

\smallskip
 \noindent
{\bf  Case~\ref{v-0-ind-b}.2.2:} $Q=C$ and $R\ne D$.

Since $w^-\sim v, v^+$, $w\not\sim v, v^+$ by Claim~\ref{nonadjacency}.
Thus, $x^+\sim v$ or $x^+\sim v^+$.
By reversing the orientation of $R$ if necessary, we assume $x^+\sim v$.
Then
$$
xz^+\oD z w\oC x \quad \mbox{and} \quad  x^+v \iR v^+w^- \iC x^+
$$
are two cycles which together contain all the  vertices in $V(C)\cup V(D)\cup V(R)$.

\qed


Let  $x\in V(C)$.  If
there exists
a cycle $D\in \sF-\{C\}$ such that $V_D(x)=B\cap V(D)$,
then we  say that $x$ is {\it bad}
w.r.t. $D$.  
Define
$$
V_{bad}=\{x\in V(C)\,|\, x \,\, \mbox{is  a bad  or A-type vertex on $C$}\}.
$$

\begin{CLA}\label{bad-vx-propty}
The vertex set $V_{bad}$ contains no two consecutive vertices on $C$.
Moreover, no other vertex in $V(C)-V_{bad}$ is adjacent to
any B-type vertex on any cycles other than $C$.
\end{CLA}

\pf Each vertex in $V_{bad}$ is adjacent to some B-type vertex
on cycles other than $C$ by the definition.
Let  $v\in V_{bad}$.
Then  by Claim~\ref{v-0-ind-b},  $v^+$ is not
adjacent to
any B-type vertex on any cycles other than $C$.
Hence, for any vertex $w\in V(C)$, $w$ or $w^+$
does not belong to $V_{bad}$. 
Thus, $V_{bad}$ contains no two consecutive vertices on $C$.


To proof the second part of the statement, assume that $v\in V(C)$
is a vertex adjacent to some B-type vertex on a cycle $D\in \sF-\{C\}$.
Since vertices in $(A\cap V(C))^+$ are not adjacent to any B-type vertices on cycles other than $C$,
$v^-$ is a B-type vertex. If $v$ is also of B-type, then
by Claim~\ref{full-b-neighbors}, $V_D(v)=B\cap V(D)$. So $v\in V_{bad}$
by the definition of $V_{bad}$. If $v$ is of A-type, then $v\in V_{bad}$
again by the definition of $V_{bad}$.
\qed

\begin{CLA}\label{a-type-wrt-c}
Let  $xy\in E(C)$ be a B-type edge. For any cycle $D\in \sF-\{C\}$,
if $V_D(xy)=B\cap V(D)$, then for any $z\in A\cap V(D)$, $z$ is of A-type w.r.t. only the cycle $C$.
\end{CLA}

\pf As $xy\in E(C)$ is of B-type, for each cycle $Q\in \sF-\{C\}$,
vertices on $Q$ are alternating between $I_{xy}$ and $V(G)-V(C)-I_{xy}$,
by $(2)$ of Claim~\ref{I-xy-alter-size-inde}.
As $I_{xy}$ is an independent set in $G$ by $(1)$ of Claim~\ref{I-xy-alter-size-inde},
and $A\cap V(D)\subseteq I_{xy}$, for any $z\in A\cap V(D)$, it is not possible for
$z$ to be adjacent to two consecutive vertices on any cycle $Q\in \sF-\{C,D\}$.
Thus, $z$ is of A-type w.r.t. only the cycle $C$.
\qed

For each vertex $x\in V_{bad}$, we define
$$
U_{x}^0=\{x^+\} \quad \mbox{and}\quad U_{x}^1=\{u\,|\,u^+\in V_C(U_{x}^0)-V_{bad}\}-U_{x}^0.
$$
For each  vertex $x_1\in U_{x}^1$, define the path
$$
P_{[x_1,x]}=x_1\iC x^+ x_1^+ \oC x
$$
to be the directed path started at $x_1$ and ended at $x$.

Start now on, if $v$ is a vertex on a directed path and
$v$ is not the end of the path, we denote by \textbf{$v^{\dag}$ }
the successor of $v$ on this path. This notation $v^{\dag}$
will be only used in the following occasion.


It is easy to see that for any  $x_2 \in V(P_{[x_1,x]})$
such that $x_2^{\dag} \in V_{P_{[x_1, x]}}(x_1)$,
$P_{[x_2, x]}=x_2\iP_{[x,x_1]}x_1x_2^{\dag}\oP_{[x,x_1]}x$ is a directed
path starting at $x_2$ and ending at $x$. Furthermore,  $P_{[x_2, x]}$ contains all the
vertices of $C$. In general, for $i\ge 2$ we define
\begin{eqnarray*}
  U_{x}^i &=& \{u\,|\, u^{\dag}\sim v,\,\,\mbox{for any $v\in U_{x}^{i-1}$, and any
  $u^{\dag}\in V(P_{[v,x]})-V_{bad}$}\}-\bigcup_{j=0}^{i-1}U_{x}^j\\
  U_{x}^{\infty} &=& \bigcup_{i=0}^{\infty}U_{x}^i.
\end{eqnarray*}

\begin{CLA}\label{zig-path-on-c}
Let $x\in V_{bad}$ and let $U_{x}^i$
be defined as above. Let $D\in \sF-\{C\}$
such that $x$ is bad or of A-type w.r.t. $D$, and
let $u\in B\cap V(D)$ such that $x\sim u$.
Then each of the following holds.
\begin{itemize}
  \item [$(1)$]$U_{x}^{\infty} \subseteq \Nb_C(uu^+)$, i.e., for any $v\in U_{x}^{\infty}$,
  $v\not\sim u, u^+$.
  \item [$(2)$] For any $y\in V(C)-V_{bad}$ such that $y$
  is adjacent to some vertex in $ U_{x}^{\infty}$, $y\sim u^+$.
  \item [$(3)$] If $x$ is  bad and $v\in U_{x}^{\infty}$, then $V_D(v)=\emptyset$.
   \item [$(4)$] If $x$ is bad and $y\in V(C)-V_{bad}$ such that $y$
  is adjacent to some vertex in $ U_{x}^{\infty}$,
  then $V_D(y)=A\cap V(D)$.
\end{itemize}
\end{CLA}

\pf We first prove $(1)$ and $(2)$ simultaneously by applying induction on
$i$.  For $i=0$, $U_{x}^{0}=\{x^+\}$.
As $x\sim u$, we have that $x^+\not\sim u^+$ by Claim~\ref{nonadjacency}.
Furthermore, as $u$ is a B-type vertex, $u\not\sim x^+$.
Hence, $x^+\in \Nb_C(uu^+)$. For any
$y\in V(C)-V_{bad}$ such that $y\sim x^+$,
since $x^+\in \Nb_C(uu^+)$, $y$ has to be  adjacent to
at least one of $u, u^+$ by the $2K_2$-freeness.
As $y\in V(C)-V_{bad}$, $y\sim u^+$ by the second part of
Claim~\ref{bad-vx-propty}. Assume now
that both $(1)$ and $(2)$ are true for $i-1$ with
$i\ge 1$. Let $v\in U_{x}^i$. By the definition
of $U_{x}^i$, there exists $w\in U_{x}^{i-1}$
such that $v\in V(P_{[w,x]})$ and $w\sim v^{\dag}$,
where $v^{\dag}\in V(C)-V_{bad}$ is the successor
of $v$ on the directed path $P_{[w,x]}$.
By the induction hypothesis, $v^{\dag}\sim u^+$.
Also, by the induction hypothesis, $U_{x}^j\subseteq \Nb_C(uu^+)$
for any $j\le i-1$. Thus, $v^{\dag}\not\in \bigcup_{j=0}^{i-1}U_{x}^j$.
Furthermore,  $v\not\in \bigcup_{j=0}^{i-1}U_{x}^j$
as $U_{x}^i$ is disjoint with $\bigcup_{j=0}^{i-1}U_{x}^j$
by its definition. Since any edge on
$P_{[w,x]}$ which is not an edge of $C$ has one
endvertex in $\bigcup_{j=0}^{i-1}U_{x}^j$, $vv^{\dag}$
is an edge on $C$. Thus, as $v^{\dag}\sim u^+$, $v\not\sim u$
by Claim~\ref{nonadjacency}.
Furthermore, $v\not\sim u^+$. For otherwise, if $v\sim u^+$,
then as $x\sim u$, and $P_{[v,x]}$ is
a spanning path of $C$, we get a cycle $v\oP_{[v,x]}xu\iD u^+ v$,
which combines $C$ and $D$ into a single cycle.
Thus, $v\in \Nb_C(uu^+)$.
For any
$y\in V(C)-V_{bad}$ such that $y\sim v$,
since $v\in \Nb_C(uu^+)$, $y$ has to be adjacent to
at least one of $u, u^+$ by the $2K_2$-freeness.
As $y\in V(C)-V_{bad}$, $y\sim u^+$ by the second part of
Claim~\ref{bad-vx-propty}.

For the statements $(3)$ and $(4)$, we see
that immediately by noticing that
the cycle $D$ is AB-alternating and
$x$ is adjacent to all the B-type vertices
on $D$ if $x$ is bad.
\qed


Define
$$
U^{\infty}=\bigcup_{x\in V_{bad}} U_{x}^{\infty}
$$
Let  $v\in U_x^{\infty}$ and $D\in \sF-\{C\}$ such that
$x$ is bad or of A-type w.r.t. $D$. Then $v$
is called {\it co-absorbable} w.r.t. $C$ and $D$ if there exists a cycle $R$
containing all the vertices in $V(C)\cup V(D)-\{v\}$.

\begin{CLA}\label{co-absorbable}
Each vertex  $v\in U_x^{\infty}$ is co-absorbable w.r.t. $C$
and a cycle $D\in \sF-\{C\}$ such that $x$ is bad or of A-type w.r.t. $D$.
\end{CLA}

\pf  If $v\in U_{x}^{0}$, then $v=x^+$.
Let $u\in B\cap V(D)$ such that $x\sim u$, and such that  $x\sim u, u^+$
if $x$ is of  A-type  w.r.t. $D$.
Then $x^+\not\sim  u^+$ by Claim~\ref{zig-path-on-c}. Furthermore,
$x^+\not\sim u$ as $u\in B\cap V(D)$.
Thus, $x^{++}\sim u$ or $x^{++}\sim u^+$. Since $D$ is AB-alternating, $u^+$
is of A-type. By Claim~\ref{a-type-wrt-c}, $u^+$
is of A-type w.r.t. only $C$.
Let $ww^+\in E(C)$ such that $u^+\sim w,w^+$.
If $x$ is bad w.r.t. $D$, then $x\sim u, u^{++}$.
And if $x\sim u, u^{++}$, then $u^+\not\sim x, x^+$ by
Claim~\ref{nonadjacency}.
Thus, $\{x, x^+\}\cap \{w,w^+\}=\emptyset$ if  $x$ is bad w.r.t. $D$.
%
Then
$$
\left\{
  \begin{array}{ll}
    xu\iD u^+ x^{++} \oC x, & \hbox{if $x^{++}\sim u^+$;}\\
     xu^+\oD u x^{++} \oC x, & \hbox{if $x^{++}\sim u$  and $x$ is of  A-type;}\\
    xu^{++}\oD u x^{++} \oC wu^+w^+\oC x, & \hbox{if $x^{++}\sim u$ and $x$ is bad.}
  \end{array}
\right.
$$
is a cycle containing all the vertices in $V(C)\cup V(D)-\{x^+\}$.

We additionally show that $x^-$ is
co-absorbable w.r.t. $C$
and $D$. (We will need this in the argument when $i\ge 1$.)
Repeat the same argument for $x^{--}$, we then have
$$
\left\{
  \begin{array}{ll}
    xu\iD u^+ x^{--} \iC x, & \hbox{if $x^{--}\sim u^+$;}\\
     xu^+\oD u x^{--} \iC x, & \hbox{if $x^{--}\sim u$  and $x$ is of  A-type;}\\
    xu^{++}\oD u x^{--} \iC w^+u^+w\iC x, & \hbox{if $x^{--}\sim u$ and $x$ is bad.}
  \end{array}
\right.
$$
is a cycle containing all the vertices in $V(C)\cup V(D)-\{x^-\}$.

Assume now that $v\in U_{x}^i$ for $i\ge 1$.
By the definition of $U_{x}^i$, we know there exists
a spanning path $P_{[v,x]}$ of $C$ with endvertices $v$ and
$x$.
By Claim~\ref{zig-path-on-c}, $v\not\sim u, u^+$.
Let $y$ be the neighbor of $v$ on $P_{[v,x]}$.
As $vy$ is an edge, and $v\not\sim u, u^+$, $y\sim u$
or $y\sim u^+$.  Since
$U_{x}^j\subseteq \Nb_C(uu^+)$
for any $j\le i-1$, we have that $y\not\in \bigcup_{j=0}^{i-1}U_{x}^j$.
Furthermore,  $v\not\in \bigcup_{j=0}^{i-1}U_{x}^j$
as $U_{x}^i$ is disjoint with $\bigcup_{j=0}^{i-1}U_{x}^j$
by its definition. Thus, $vy$ is an
edge on $C$,
since any edge on
$P_{[v,x]}$ which is not an edge of $C$ has one
endvertex in $\bigcup_{j=0}^{i-1}U_{x}^j$.
We may assume
that $y\not\in V_{bad}$, as  both the predecessor and  successor
of  a bad vertex on $C$
is co-absorbable by the argument for $i=0$ case.
Thus $y\sim u^+$ by (2) of Claim~\ref{zig-path-on-c}.
Then
$y\iP_{[v,x]} xu \iD u^+ y$ is the desired cycle.
\qed

\begin{CLA}\label{numner-of-neighbors-for-co-absorbable-vertex}
We may assume that each vertex in $U^{\infty}$ has less than $(|V(G)|-1)/3$
neighbors in $G$.
\end{CLA}

\pf Suppose on the contrary that there exists $v\in U^{\infty}$
so that $|N_G(v)|\ge (|V(G)|-1)/3$. By Claim~\ref{co-absorbable},
we see that $v$  is co-absorbable w.r.t. $C$ and some cycle $D\in \sF-\{C\}$.
By standard arguments for longest cycles,  we know that
$v$ has no two neighbors which are consecutive on any cycle
$Q\in\sF-\{C,D\}$ and on the cycle which is the combination of $C-v$ and $D$;
and also that $(N_G(v))^+$,
the set of the successors  of  neighbors of $v$
from the cycle which is the combination of  $C-v$ and $D$
and cycles in $\sF-\{C,D\}$, is an independent set in $G$.
Let $S=V(G)-(N_G(v))^+-\{v\}$. Then
$c(G-S)=|(N_G(v))^{+}\cup \{v\}|\ge (|V(G)|-1)/3+1>\frac{|V(G)|}{3}$.
So $\tau(G)\le \frac{|S|}{c(G-S)}< 2$. This achieves a
contradiction.
\qed

\begin{CLA}\label{zig-indep}
Each of the following holds.
\begin{itemize}
  \item [$(1)$]The set $U^{\infty}$  is an independent set in $G$.
  \item [$(2)$]$V_C(U^{\infty})\cap U^{\infty}=\emptyset$.
  \item [$(3)$]$ U^{\infty}\cup A^+$ is an independent set in $G$.
\end{itemize}
\end{CLA}

\pf To prove (1),
assume that there exist $u, v\in U^{\infty}$
such that $uv\in E(G)$.
By Claim~\ref{numner-of-neighbors-for-co-absorbable-vertex},
$u$ and $v$ in total have at most
$2(|V(G)|-1)/3$ neighbors in $G$.
As $uv$ is an edge, and $G$ is $2K_2$-free,
the set of non-neighbors of $u$ and $v$
in $G$ forms an independent set in $G$.
Let $S=N_G(u)\cup N_G(v)-\{u\}$.
Then $c(G-S)=|V(G)-S-\{u\}|>|V(G)|/3$.
So $\tau(G)<2$.
Again, we achieve  a contradiction to
the assumption that $\tau(G)\ge 3$.
As $U^{\infty}$  is an independent set in $G$, we have
$V_C(U^{\infty})\cap U^{\infty}=\emptyset$.
Since each bad vertex $x$ is adjacent to its successor $x^+$, and 
 $x^+\in U_x^0\subseteq U^{\infty}$,
we have that $V_{bad}\subseteq V_C(U^{\infty})$.
Thus,
 no vertex in $U^{\infty}$
is adjacent to any B-type vertex on cycles other than $C$.
Since $(A\cap V(C))^+\subseteq U^{\infty}$, we know that
$U^{\infty}\cup A^+$ is an independent set in $G$.
\qed


\begin{CLA}\label{S-comp-correspondence}
For any vertex $y\in V_C(U^{\infty})$, there exists $v\in U^{\infty}$
such that $vy\in E(C)$.
\end{CLA}

\pf Assume that $y\in V_C(U_x^{\infty})$ for
some $x\in V_{bad}$. The Claim trivially holds
if $y\in V_C(U_x^{0})$. So assume that $i\ge 1$
and let $y\in V_C(U_x^i)-V_C(\bigcup_{j=0}^{i-1} U_x^j)$.
By the definition of $U_x^{i}$, we know that
there exists $w\in U_x^{i}$, and a spanning path $P_{[w,x]}$
of $C$ with endvertices as $w$ and $x$
such that $y$  is a neighbor of $w$
on $P_{[w,x]}$. Since $V_C(U^{\infty})\cap U^{\infty}=\emptyset$ by (2)
of Claim~\ref{zig-indep}, $y\not\in U_x^{\infty}$.
By the assumption that $y\in V_C(U_x^i)-V_C(\bigcup_{j=0}^{i-1} U_x^j)$,
we know that the predecessor $v$ of $y$
on $P_{[w,x]}$ satisfies that
$v\not\in \bigcup_{j=0}^{i-1} U_x^{j}$.
As any edge of $P_{[w,x]}$ which is not
an edge of $C$ has one end contained in $\bigcup_{j=0}^{i-1} U_x^{j}$,
we then know that $vy\in E(C)$.
%
\qed

\begin{CLA}\label{S-comp-correspondence2}
$|V_C(U^{\infty})|\le 2|U^{\infty}|$.
\end{CLA}

\pf Since $U^{\infty}$ is an independent set in $G$ by
Claim~\ref{zig-indep}, $|N_C(U^{\infty})|\le 2|U^{\infty}|$.
Let $y\in V_C(U^{\infty})$ be any vertex. By Claim~\ref{S-comp-correspondence},
there exists $v\in U^{\infty}$
such that $vy\in E(C)$.
Thus, $V_C(U^{\infty})\subseteq N_C(U^{\infty})$.
So $|V_C(U^{\infty})|\le |N_C(U^{\infty})|\le 2|U^{\infty}|$.
\qed

Let
$$
S=A\cup V_C(U^{\infty}).
$$
We claim that each vertex in $A^+\cup U^{\infty}$
is an isolated vertex in $G-S$.
This is because $A^+\cup U^{\infty}$ is an independent
set in $G$, and all the possible neighbors of
vertices in $A^+\cup U^{\infty}$ in $G$
are contained in $S$.
Note also that $|V_C(U^{\infty})|\le 2|U^{\infty}|$ by Claim~\ref{S-comp-correspondence2},
and $|S\cap (V(G)-V(C))|=|V(G)-V(C)-S|=|V(G)-V(C)|/2$
as we assume that all cycles in $\sF-\{C\}$ are AB-alternating.
Since $A\cap V(C)\subseteq V_{bad}$ by the definition of $V_{bad}$,
and $V_{bad}\subseteq V_C(U^{\infty})$, we have that $A\cap V(C)\subseteq V_C(U^{\infty})$.
Thus, $S=A\cup V_C(U^{\infty})=V_C(U^{\infty})\cup (A\cap (V(G)-V(C)))$ and 
thus $|S|=|V_C(U^{\infty})|+|V(G)-V(C)|/2$.
Hence
\begin{eqnarray*}
  \tau(G) &\le & \frac{|S|}{c(G-S)} \\
   &\le & \frac{|V_C(U^{\infty})|+|V(G)-V(C)|/2}{|U^{\infty}|+|V(G)-V(C)|/2} \\
   &\le & \frac{2|U^{\infty}|+|V(G)-V(C)|/2}{|U^{\infty}|+|V(G)-V(C)|/2}<2,
\end{eqnarray*}
showing a contradiction.
The proof of Theorem~\ref{main} is now complete.
\qqed

\textbf{{\noindent \large Acknowledgements}}

The author is extremely  grateful to
Professor Mark Ellingham for his careful
comments and suggestions in improving the
proofs and the writing of this paper.

\bibliographystyle{plain}
\bibliography{SSL-BIB}

\end{document}